\documentclass[11pt,a4paper]{article}
\usepackage{pifont}
\usepackage{CJK,epsfig,amsfonts,amsgen,amsmath,amstext,amsbsy,amsopn,amsthm
}
\usepackage{mathrsfs}
\usepackage{float}
\usepackage{epsfig,latexsym,amsfonts,amsmath,amssymb}
\usepackage{color}
\usepackage{calc}
\usepackage{enumerate}

\usepackage{comment}
\usepackage{geometry}
\usepackage{calc}
\usepackage{multirow}

\usepackage{tikz}
\usetikzlibrary{shapes.geometric, arrows, positioning}
\usetikzlibrary{decorations.markings}
\tikzstyle{vertex}=[circle, draw, inner sep=0pt, minimum size=6pt]
\tikzstyle{node} = [rectangle, minimum width=1.5cm, minimum height=0.7cm, text centered, draw=black, fill=blue!20, font=\small]
\tikzstyle{arrow} = [thick,->,>=stealth]

\newtheorem{theorem}{Theorem}
\newtheorem{lemma}[theorem]{Lemma}

\newtheorem{corollary}[theorem]{Corollary}
\newtheorem{remark}[theorem]{Remark}

\newcommand{\Fq}{\mathbb F_q}
\newcommand{\gbinom}[2]{\genfrac{[}{]}{0pt}{}{#1}{#2}}

\newcommand{\one}{\mathbf 1}

\usepackage{indentfirst} 
\newcounter{Hcase}

\theoremstyle{plain}

\date{\today}
\begin{document}

\title{A $d$-Degree Erd\H{o}s--Ko--Rado Theorem for Cross-Intersecting Families in Finite Vector Spaces}

\author{\small Zhen Jia${\thanks{Research partially supported by the National Key R\&D Program of China (No. 2025YFA1017700). E-mail address: 12231271@mail.sustech.edu.cn}}^{~1}$,~Qing Xiang${\thanks{Research partially supported by the National Key R\&D Program of China (No. 2025YFA1017700) and the National Natural Science Foundation of China Grant (No. 12131011). E-mail address: xiangq@sustech.edu.cn}}^{~1}$, and Jimeng Xiao${\thanks{Research partially supported by the National Natural Science Foundation of China Grant (No. 12501460). E-mail address: xiaojm@sustech.edu.cn}}^{~1}$
 \\[2mm]
\small ${}^1$  Department of Mathematics\\[-0.8ex]
\small Southern University of Science and Technology\\[-0.8ex]
\small Shenzhen, Guangdong 518055, China}

\baselineskip 20pt

\date{}
\maketitle
 \vspace{4mm}

\begin{abstract}
Fix integers $1\le k\le n$. Let $q$ be a prime power, let $V$ be an $n$-dimensional vector space over $\mathbb F_q$, and let
$\gbinom{V}{k}$ denote the family of all $k$-dimensional subspaces of $V$.
For a family $\mathcal F$ of $k$-dimensional subspaces of $V$ and $1\leq d<k$, define the minimum $d$-degree of $\mathcal F$ by
$ \delta_d(\mathcal F)
    =
    \min_{S\in \gbinom{V}{d}}|\{F\in \mathcal F:S\le F\}|$.
We prove that if
$1\le d<k$, $n\ge 2k$, and
$\mathcal F,\mathcal G\subseteq \gbinom{V}{k}$ are cross-intersecting, then
\[
    \delta_d(\mathcal F)\delta_d(\mathcal G)
    \le
    \gbinom{n-d-1}{k-d-1}^2.
\]
The result holds uniformly for all \(1\le d<k\), with the sharp threshold \(n\ge 2k\) independent of \(d\).

\vskip 0.1in \noindent%
\textbf{Keywords}: Erd\H{o}s--Ko--Rado theorem, cross-intersecting families, vector space, $d$-degree, the eigenvalue method.
 \\[7pt]
{\sl MSC:}\ \ 05D05
\end{abstract}

\section{Introduction}
Fix integers $1\le k\le n$. Let $\mathcal{F}$ be a family of $k$-subsets of $[n]:=\{1,2,\ldots,n\}$. We say that $\mathcal{F}$ is \emph{intersecting} if any two $k$-subsets in $\mathcal{F}$ have nonempty intersection. The celebrated Erd\H{o}s--Ko--Rado theorem~\cite{erdos1961} states that when $n\geq 2k$, every intersecting family $\mathcal{F}$ of $k$-subsets of $[n]$
must satisfy $|\mathcal{F}|\leq \binom{n-1}{k-1}$. Equality can be reached if one chooses all $k$-subsets of $[n]$ containing a fixed element. Furthermore, if $n>2k$, then these are the only extremal families.
The Erd\H{o}s--Ko--Rado theorem is a cornerstone of extremal combinatorics and has many interesting
applications and generalizations. For a comprehensive survey,
see \cite{frankTokushigel2016}.

In \cite{huang2017}, Huang and Zhao initiated the study of intersecting families $\mathcal{F}\subseteq \binom{[n]}{k}$ from the perspective of minimum degree. They established the following degree version of the Erd\H{o}s--Ko--Rado theorem.
\begin{theorem}[\cite{huang2017}]\label{huang-1-degree}
    Let $k,n$ be positive integers with $n\geq 2k+1$. Suppose $\mathcal{F}\subseteq \binom{[n]}{k}$ is an intersecting family. Define $\delta_1(\mathcal{F})=min_{i\in [n]}|\{F\in \mathcal{F}:i\in F\}|$. Then
    \[
    \delta_1(\mathcal{F})\leq \binom{n-2}{k-2}.
    \]
    Moreover, equality holds if and only if $\mathcal{F}=\left\{F\in\binom{[n]}{k}:a\in F\right\}$ for some $a\in [n]$.
\end{theorem}

We remark that when $n\ge 2k+1$, Theorem~1 implies the EKR theorem. The authors of [10]  also proved the following cross-intersecting minimum-degree version, which in particular generalizes Theorem \ref{huang-1-degree}.

\begin{theorem}[\cite{huang2017}]
    For $n \geq 2k + 1$, two cross-intersecting families $\mathcal{F}$ and $\mathcal{G}$ of $k$-subsets of $[n]$ satisfy $\delta_{1}(\mathcal{F})\delta_{1}(\mathcal{G})\leq \binom{n-2}{k-2}^2$.
\end{theorem}

For a family $\mathcal F$ of $k$-subsets of $[n]$, and $1\leq d< k$, define 
\[
\delta_{d}(\mathcal{F})=\min_{|S|=d,S\subseteq [n]}|\left\{F\in \mathcal{F}:S\subseteq F\right\}|.
\]

For $d\geq 2$, Kupavskii \cite{Kupavskii2019} showed that $\delta_{d}(\mathcal{F})\leq \binom{n-d-1}{k-d-1}$ for $n\geq 2k+\frac{3d}{1-\frac{d}{k}}$. Later, Huang and Zhang \cite{huang2026} used techniques from spectral graph theory to improve the sufficient range to $n\geq 2k+2d-3$. In \cite{Han2026}, Gan, Han and Im proved that the same bound holds for $n\geq 2k+\sqrt{2k}+O(1)$ when $d=k-1$, and for $n\geq 2k+7k^{2/3}+O(k^{1/3})$ when $d=k-2$.

Problems in extremal set theory have natural vector space analogues, or $q$-analogues. The study of the 
vector space analogue of the EKR theorem was initiated by
Hsieh~\cite{Hsieh1975}; subsequently, Frankl and Wilson~\cite{FranklWilson1986} used spectral methods to prove the EKR theorem for $t$-intersecting families of subspaces. For further results, see
\cite{caolv2023,caolv2022}.

We use the following notation throughout this paper. Let $q$ be a prime power and let $V$ be an $n$-dimensional vector space over $\Fq$.  For $0\le k\le n$, denote by
$\gbinom{V}{k}_q$ the collection of all $k$-dimensional subspaces of $V$. Let $\mathcal{F}\subseteq\gbinom{V}{k}_q$ be a family of subspaces. We say that $\mathcal{F}$ is \emph{intersecting} if for all $F,F'\in\mathcal{F}$, we have $\dim(F\cap F')\geq 1$. Similarly, two families $\mathcal{F},{\mathcal{G}}\subseteq\gbinom{V}{k}_q$ are \emph{cross-intersecting} if $\dim(F\cap G)\geq 1$, for all $F\in\mathcal{F}$ and $G\in \mathcal{G}$.

For integers $0\le r\le m$, the \emph{$q$-binomial coefficient}, or \emph{Gaussian binomial coefficient} $\gbinom{m}{r}_q$ is defined by
\[
    \gbinom{m}{r}_q
    =
    \prod_{i=0}^{r-1}
    \frac{q^{m-i}-1}{q^{r-i}-1}.
\]
We use the convention $\gbinom{m}{0}_q=1$ and $\gbinom{m}{r}_q=0$ if $r<0$ or $r>m$. To simplify notation, we shall omit the subscript $q$. 

Frankl and Tokushige \cite{frankTokushigel2017} proved the vector space analogue of the minimum degree EKR theorem. Recently, Shan and Zhou \cite{ShanZhou2025} established the vector space version of the minimum $d$-degree EKR theorem. For a family $\mathcal F\subseteq \gbinom{V}{k}$ and an integer $d$ with $1\le d<k$, define the minimum $d$-degree of
$\mathcal F$ by
\[
\delta_d(\mathcal F)
    =
    \min_{S\in \gbinom{V}{d}}|\{F\in \mathcal F:S\le F\}|.
\]
\begin{theorem}[\cite{frankTokushigel2017}]
    Let $n,k$ be positive integers with $n\geq 2k\geq 4$. Suppose $\mathcal F\subseteq \gbinom{V}{k}$ is an intersecting family. Then
    \[
    \delta_{1}(\mathcal{F})\leq \gbinom{n-2}{k-2}.
    \]
\end{theorem}
\begin{theorem}[\cite{ShanZhou2025}]
    Let $n,k,d$ be positive integers with $k>d\geq 2$ and $n\geq 2k+1$. Suppose $\mathcal F\subseteq \gbinom{V}{k}$ is an intersecting family. Then
    \[
    \delta_{d}(\mathcal{F})\leq \gbinom{n-d-1}{k-d-1}.
    \]
\end{theorem}
Shan and Zhou~\cite{ShanZhou2025} asked whether spectral methods could be used to establish an analogous cross-intersecting minimum $d$-degree bound for vector spaces. 
In this paper, we answer this question in the affirmation.

\begin{theorem}\label{thm:main}
Let $1\le d<k$ and $n\ge 2k$.  Suppose that
$\mathcal F,\mathcal G\subseteq \gbinom{V}{k}$
are cross-intersecting.  Then
\[
    \delta_d(\mathcal F)\delta_d(\mathcal G)
    \le
    \gbinom{n-d-1}{k-d-1}^2.
\]
\end{theorem}
\begin{remark}
The upper bound in Theorem~\ref{thm:main} is attained by
\[
\mathcal F=\mathcal G=
\left\{
K\in\gbinom{V}{k}:L\le K
\right\},
\]
where \(L\le V\) is a fixed one-dimensional subspace of $V$.

The condition \(n\ge2k\) is best possible and does not depend on the value of $d$. Indeed, if
\(k<n<2k\), then any two \(k\)-dimensional subspaces of \(V\) have
a nontrivial intersection, since
$\dim(F\cap G)\ge 2k-n\ge1$, for any $F,G\in \gbinom{V}{k}$. 
Hence the families
\(\mathcal F,\mathcal G\subseteq\gbinom{V}{k}\)
are cross-intersecting. In this range, the trivial upper bound becomes sharp:
\[
\delta_d(\mathcal F)\delta_d(\mathcal G)
\le
\gbinom{n-d}{k-d}^{2},
\]
and equality is attained by taking $\mathcal F=\mathcal G=\gbinom{V}{k}$.

We also emphasize that Theorem~\ref{thm:main} treats all
\(1\le d<k\) in a unified way. In the proof given in Section 2, the exact form of one auxiliary spectral estimate differs slightly between the cases \(d=1\) and
\(d\ge2\); however, only a common upper bound is needed for the
final argument. Therefore, the theorem does not need to be divided
into separate cases according to \(d\).
\end{remark}

\noindent\textbf{Proof Overview and Paper Organization.}
Our proof relies on a sharp two-sided estimate for $\Omega$, a weighted count of pairs of disjoint $d$-dimensional subspaces formally defined in Section~\ref{sec:preliminaries}.
This approach builds on the spectral framework of Huang and Zhao~\cite{huang2017}, the degree-counting techniques of Huang and Zhang~\cite{huang2026}, and the vector-space spectral calculations of Frankl-Tokushige~\cite{frankTokushigel2017} and Shan-Zhou~\cite{ShanZhou2025}.
To reach the optimal threshold $n \ge 2k$, we develop a sharper control over the spectral error terms.
The primary innovation is a refined estimate for $R_{\mathcal F}R_{\mathcal G}$, which measures the contribution from the eigenspaces beyond $U_0 \oplus U_1$.
After eliminating the $U_1$-term, this yields the upper bound on $\Omega$ in Corollary~\ref{lem:spectral-upper}.
On the combinatorial side, the
balancing factor $c=\sqrt{{\delta_d(\mathcal G)}/{\delta_d(\mathcal F)}}$
symmetrizes the two families and turns the assumed lower bound on
\(\delta_d(\mathcal F)\delta_d(\mathcal G)\) into the lower estimate for
\(\Omega\) in Corollary~\ref{lowerbound-omerga}. 
The alignment of these two spectral and combinatorial refinements enables a unified treatment of all $1 \le d < k$ and successfully covers the boundary case $n = 2k$.
The remainder of the paper is structured as follows: Section~\ref{sec:preliminaries} develops the necessary spectral tools, while Section~\ref{section-proof-mainth} establishes the two-sided bounds on $\Omega$ and completes the proof of Theorem~\ref{thm:main} via a quadratic factorization argument.


\section{Preliminaries}\label{sec:preliminaries}

We will need  the following elementary identities and counting formulas.
\[
    \frac{\gbinom{m-1}{r-1}}{\gbinom{m}{r}}
    =
    \frac{q^r-1}{q^m-1},
    \qquad
    \frac{\gbinom{m-1}{r}}{\gbinom{m}{r}}
    =
    \frac{q^{m-r}-1}{q^m-1}.
\]
Equivalently,
\[
    \gbinom{m}{r}
    =
    \frac{q^m-1}{q^r-1}\gbinom{m-1}{r-1}
    =
    \frac{q^m-1}{q^{m-r}-1}\gbinom{m-1}{r}.
\]
We also have
\[
    \gbinom{n}{d}\gbinom{n-d}{k-d}
    =
    \gbinom{n}{k}\gbinom{k}{d}.
\]

Recall that $V$ is an $n$-dimensional vector space over $\mathbb F_q$.

\begin{lemma}\cite[Lemma~9.3.2]{GodsilMeagher2016}\label{lem:disjoint-count}
Let $S\le V$ be a fixed $d$-dimensional subspace and $n\geq 2d$.  Then the number of
$d$-dimensional subspaces $T\le V$ satisfying
$S\cap T=\{0\}$
is $q^{d^2}\gbinom{n-d}{d}$.
\end{lemma}

Let $\Gamma=\Gamma_q(n,k)$ be the $q$-Kneser graph on $\gbinom{V}{k}$. 
Thus the vertex set is $\gbinom{V}{k}$ and two vertices $F,F'\in \gbinom{V}{k}$ are adjacent if
$F\cap F'=\{0\}$.
Let $A_{\Gamma}$ be the adjacency matrix of $\Gamma$, and set
\begin{equation*}
    A:=q^{-k^2}A_{\Gamma}
\end{equation*}
The following spectral decomposition is standard.

\begin{lemma}\cite[Corollary~9.7.1]{GodsilMeagher2016}\label{lem:q-kneser-spectrum}
The vector space $\mathbb R^{\gbinom{n}{k}}$ has an orthogonal decomposition
\[
    \mathbb R^{\gbinom{n}{k}}
    =
    U_0\oplus U_1\oplus\cdots\oplus U_k
\]
where $U_j$ is the eigenspace of $A$ of the eigenvalue
\[
    \lambda_j
    =
    (-1)^j q^{\binom j2-kj}
    \gbinom{n-k-j}{k-j},
\]
with multiplicity $\gbinom{n}{j}-\gbinom{n}{j-1}$, for any $0\le j\le k$.
Moreover, $U_0$ is spanned by the all-ones vector $\one$.
\end{lemma}

In particular, 
$\lambda_0=\gbinom{n-k}{k}>0, \lambda_1=-q^{-k}\gbinom{n-k-1}{k-1}<0,\text{and } \lambda_2=q^{1-2k}\gbinom{n-k-2}{k-2}>0$ when $k\ge 2$. We shall need the following result about the eigenvalues mentioned above.

\begin{lemma}\label{lem:lambda-comparison}
Assume $n\ge 2k$ and $k\ge 2$.  Then
the sequence $(|\lambda_j|)_{j\geq 1}$ is strictly decreasing.
Moreover,
\[
    r:=\frac{\lambda_2}{|\lambda_1|}
    =
    q^{1-k}
    \frac{q^{k-1}-1}{q^{n-k-1}-1}
    \leq\frac{1}{2}.
\]
\end{lemma}

\begin{proof}
For $1\le j<k$, we compute
\begin{equation*}
    \frac{|\lambda_{j+1}|}{|\lambda_j|}=q^{j-k}
    \frac{\gbinom{n-k-j-1}{k-j-1}}
         {\gbinom{n-k-j}{k-j}}=q^{j-k}
    \frac{q^{k-j}-1}{q^{n-k-j}-1}.
\end{equation*}
Since $n\ge 2k$, we have $q^{n-k-j}-1\ge q^{k-j}-1$, and therefore $ \frac{|\lambda_{j+1}|}{|\lambda_j|}
\le q^{j-k}< 1$.
Thus the sequence $|\lambda_j|$ is decreasing for $j\ge 1$.

Finally, since $n\ge 2k$, we have 
\begin{equation*}
 r=\frac{\lambda_2}{|\lambda_1|}=\frac{
    q^{1-2k}\gbinom{n-k-2}{k-2}
    }{q^{-k}\gbinom{n-k-1}{k-1}}=
    q^{1-k}
    \frac{q^{k-1}-1}{q^{n-k-1}-1} \le q^{1-k}\le \frac{1}{2}.
\end{equation*} 
\end{proof}

Let $\mathcal F,\mathcal G\subseteq \gbinom{V}{k}$, and
let $\chi_{\mathcal F},\chi_{\mathcal G}\in\mathbb R^{\gbinom{n}{k}}$ be their
characteristic vectors. That is, $\chi_{\mathcal F}$ and $\chi_{\mathcal G}$ are $\{0,1\}$-valued column vectors of dimension $\gbinom{n}{k}$, 
and for each $Y\in\gbinom{V}{k}$, the $Y^{\rm th}$-coordinate of $\chi_{\mathcal F}$ (resp. $\chi_{\mathcal G}$) is $1$ if $Y\in\mathcal F$ (resp.\ $Y\in\mathcal G$), and $0$ otherwise. By Lemma \ref{lem:q-kneser-spectrum}, we write
\[
    \chi_{\mathcal F}=a_0+a_1+\cdots+a_k,
    \quad
    \chi_{\mathcal G}=b_0+b_1+\cdots+b_k, \text{ where }a_j,b_j\in U_j, \;\forall j, 0\le j\le k.
\]
For $S\in\gbinom{V}{d}$, we define $d_S(\mathcal F)=|\left\{F\in\mathcal{F}:S\leq F \right\}|$ and
\[
    \Omega
    =
    \sum_{\substack{S,T\in\gbinom{V}{d}\\S\cap T=\{0\}}}
    d_S(\mathcal F)d_T(\mathcal G).
\]

\begin{lemma}\cite[Lemma~2.7]{ShanZhou2025}\label{lem:mu-expression}
Assume that $1\leq d<k$ and $n\geq 2k$. Let $\mathcal{F},\mathcal{G}\subseteq \gbinom{V}{k}$ be cross-intersecting. Then
we have
\[
    \Omega
    =
    \sum_{j=0}^{d}\mu_j a_j^{\top}b_j,
\]
where
\[
    \mu_j
    =
    (-1)^j
    q^{kj+(d-j)^2-\binom{j+1}{2}}
    \gbinom{k-j}{d-j}
    \gbinom{n-d-j}{d-j}
    \gbinom{n-d-j}{k-d}
\]
for $0\le j\le d$.  For convenience, we put $\mu_j=0$ for $j>d$.
\end{lemma}
For the convenience of the reader, we include a proof of Lemma \ref{lem:mu-expression}.
\begin{proof}
    For $0\leq i\leq j\leq k$, let $W_{i,j}$ denote the $\gbinom{n}{i}\times \gbinom{n}{j}$ matrix whose rows are indexed by the $i$-dimensional subspaces of $V$, columns are indexed by the $j$-dimensional subspaces of $V$, and 
    \[
(W_{i,j})_{S,T} =
\begin{cases}
    1,      & \text{if } S \leq T, \\
    0,       & \text{otherwise.}
\end{cases}
\]
Similarly, we define the matrices $\overline{W}_{i,j}$ of size $\gbinom{n}{i}\times\gbinom{n}{j}$ by 
    \[
(\overline{W}_{i,j})_{S,T} =
\begin{cases}
    1,      & \text{if }S \cap T=\{0\}, \\
    0,       & \text{otherwise.}
\end{cases}
\]
Let $\chi_{\mathcal F},\chi_{\mathcal G}\in\mathbb R^{\gbinom{n}{k}}$ be the characteristic vectors of $\mathcal{F}$ and $\mathcal{G}$, respectively. Write $\chi_{\mathcal F}=a_0+a_1+\cdots+a_k$ and $\chi_{\mathcal G}=b_0+b_1+\cdots+b_k$ with $a_j,b_j\in U_j$ where $0\leq j\leq k$. Then we obtain
\begin{equation*}
\Omega
    =
    \sum_{{S,T\in\gbinom{V}{d},S\cap T=\{0\}}}
    d_S(\mathcal F)d_T(\mathcal G)=\chi_{\mathcal{F}}^{\top}W_{d,k}^{\top}\overline{W}_{d,d}W_{d,k}\chi_{\mathcal{G}}.   
\end{equation*}
The calculations in the proof of
\cite[Lemma~2.7]{ShanZhou2025} shows that
$W_{d,k}^{\top}\overline{W}_{d,d}W_{d,k}$ acts on $U_j$ as multiplication by $\mu_j$ when $0\leq j\leq d$, and vanishes on $U_j$ when $j>d$. Hence, by the orthogonality of the subspaces $U_0,\ldots,U_k$, we have
\begin{equation*}
 \Omega=\sum_{i,j=0}^{k}a_i^{\top}(W_{d,k}^{\top}\overline{W}_{d,d}W_{d,k})b_j=\sum_{j=0}^{d}\mu_ja_j^{\top}b_j.  
\end{equation*}

\end{proof}

For the standard cross-Hoffman product bound, we refer the reader to \cite{Tokushige2013}. 
The following lemma contains both the Tokushige bound (\cite[Lemma~1]{Tokushige2013}) and the refined estimate for the error term; 
a proof is provided for completeness. 

\begin{lemma}\label{lem:cross-stability}
Assume that $n\geq 2k$ and $k\geq 2$.
Let $\mathcal F,\mathcal G\subseteq \gbinom{V}{k}$ be cross-intersecting families. Let $\chi_{\mathcal F},\chi_{\mathcal G}\in\mathbb R^{\gbinom{n}{k}}$ be the characteristic vectors of $\mathcal{F}$ and $\mathcal{G}$, respectively.
Write $\chi_{\mathcal F}=a_0+a_1+\cdots+a_k$ and $\chi_{\mathcal G}=b_0+b_1+\cdots+b_k$ with $a_j,b_j\in U_j$ where $0\leq j\leq k$. Define 
\[
    R_{\mathcal F}:=\left(\sum_{j=2}^k\|a_j\|^2\right)^{1/2},
    \qquad
    R_{\mathcal G}:=\left(\sum_{j=2}^k\|b_j\|^2\right)^{1/2}.
\]
Then
\[
    \sqrt{|\mathcal F| |\mathcal G|}\le \gbinom{n-1}{k-1},\text{ and } R_{\mathcal F}R_{\mathcal G}
    \le
    \frac{\sqrt{|\mathcal F| |\mathcal G|}(\gbinom{n-1}{k-1}-\sqrt{|\mathcal F| |\mathcal G|})}{\gbinom{n-1}{k-1}(1-\frac{\lambda_2}{|\lambda_1|})}.
\]
\end{lemma}

\begin{proof}
Since $\mathcal F$ and $\mathcal G$ are cross-intersecting, there are no edges of
the $q$-Kneser graph between $\mathcal F$ and $\mathcal G$. It follows that 
\[
    \chi_{\mathcal F}^{\top}A\chi_{\mathcal G}=0.
\]
Using the eigenspace decomposition in Lemma~\ref{lem:q-kneser-spectrum}, we obtain
\[
    \sum_{j=0}^k \lambda_j a_j^{\top}b_j=0.
\]

Put $X:=\sqrt{|\mathcal F| |\mathcal G|}$ and $N:=\gbinom{n}{k}$. Since $U_0$ is spanned by $\one$, we have
\[
    a_0=\frac{|\mathcal{F}|}{N}\one,
    \qquad
    b_0=\frac{|\mathcal{G}|}{N}\one,
\]
and hence $a_0^{\top}b_0=\frac{X^2}{N}$.
Thus
\[
    0
    =
    \lambda_0\frac{X^2}{N}
    +
    \sum_{j=1}^k\lambda_j a_j^{\top}b_j.
\]
By Lemma~\ref{lem:lambda-comparison},
\[
    |\lambda_j|\le |\lambda_1|
    \qquad
    \text{for }j\ge 1.
\]
Therefore, by the Cauchy-Schwarz inequality, 
\begin{equation*}
\lambda_0\frac{X^2}{N}\le
    |\lambda_1|
    \left(
    \sum_{j=1}^k\|a_j\|^2
    \right)^{1/2}
    \left(
    \sum_{j=1}^k\|b_j\|^2
    \right)^{1/2}=
    |\lambda_1|
    \sqrt{
    \left(|\mathcal{F}|-\frac{|\mathcal{F}|^2}{N}\right)
    \left(|\mathcal{G}|-\frac{|\mathcal{G}|^2}{N}\right)
    }.    
\end{equation*}
Since $\sqrt{(1-\alpha)(1-\beta)}
    \le
    1-\sqrt{\alpha\beta}$ for $0\le \alpha,\beta\le 1$,
we get $\sqrt{
    \left(|\mathcal{F}|-\frac{|\mathcal{F}|^2}{N}\right)
    \left(|\mathcal{G}|-\frac{|\mathcal{G}|^2}{N}\right)
    }
    \le
    X\left(1-\frac{X}{N}\right)$.
Since $X>0$, dividing both sides by $X$ gives $$\lambda_0 X/N \le |\lambda_1|(1-X/N);$$
that is,
$X\le
    \frac{|\lambda_1|N}{\lambda_0+|\lambda_1|}$.
Using $\lambda_0=\gbinom{n-k}{k}$ and $|\lambda_1|=q^{-k}\gbinom{n-k-1}{k-1}$, 
we obtain
\[
    X\le \gbinom{n-1}{k-1}.
\]

This is the Tokushige bound from \cite[Lemma~1]{Tokushige2013}. It remains to prove the refined estimate for
$R_{\mathcal F}R_{\mathcal G}$.
Starting again from
\[
    0
    =
    \lambda_0\frac{X^2}{N}
    +
    \lambda_1a_1^{\top}b_1
    +
    \sum_{j=2}^k\lambda_j a_j^{\top}b_j,
\]
and using $\lambda_1<0$ and $|\lambda_j|\le \lambda_2$ for all $j\ge 2$, we obtain
\[
    \lambda_0\frac{X^2}{N}
    \le
    |\lambda_1|\|a_1\|\|b_1\|
    +
    \lambda_2 R_{\mathcal F}R_{\mathcal G}.
\]

Now
\[
    \|a_1\|^2
    =
    \|\chi_{\mathcal F}\|^2-\|a_0\|^2-\sum_{j=2}^k\|a_j\|^2
    =
    |\mathcal{F}|-\frac{|\mathcal{F}|^2}{N}-R_{\mathcal F}^2,
\]
and similarly $\|b_1\|^2
    =
    |\mathcal{G}|-\frac{|\mathcal{G}|^2}{N}-R_{\mathcal G}^2$.
We use the elementary inequality
\[
    \sqrt{(A-u^2)(B-v^2)}
    \le
    \sqrt{AB}-uv
\]
which is valid for $A,B,u,v\ge 0$ with $u^2\le A$ and $v^2\le B$. 
In the above inequality, put
\[
    A=|\mathcal{F}|-\frac{|\mathcal{F}|^2}{N},
    \qquad
    B=|\mathcal{G}|-\frac{|\mathcal{G}|^2}{N},
    \qquad
    u=R_{\mathcal F},
    \qquad
    v=R_{\mathcal G},
\]
we obtain
\begin{equation*}
 \|a_1\|\|b_1\|  \le
    \sqrt{
    \left(|\mathcal{F}|-\frac{|\mathcal{F}|^2}{N}\right)
    \left(|\mathcal{G}|-\frac{|\mathcal{G}|^2}{N}\right)
    }
    -
    R_{\mathcal F}R_{\mathcal G} 
    \le
    X\left(1-\frac{X}{N}\right)
    -
    R_{\mathcal F}R_{\mathcal G}.
\end{equation*}

Consequently,
\[
\begin{aligned}
    \lambda_0\frac{X^2}{N}
    &\le
    |\lambda_1|
    X\left(1-\frac{X}{N}\right)
    -
    |\lambda_1|R_{\mathcal F}R_{\mathcal G}
    +
    \lambda_2R_{\mathcal F}R_{\mathcal G}       \\
    &=
    |\lambda_1|
    X\left(1-\frac{X}{N}\right)
    -
    (|\lambda_1|-\lambda_2)
    R_{\mathcal F}R_{\mathcal G}.
\end{aligned}
\]
Thus
\[
    (|\lambda_1|-\lambda_2)R_{\mathcal F}R_{\mathcal G}
    \le
    |\lambda_1|X
    -
    \frac{|\lambda_1|+\lambda_0}{N}X^2.
\]
Since $ \frac{|\lambda_1|N}{|\lambda_1|+\lambda_0}=\gbinom{n-1}{k-1}$
and $|\lambda_1|-\lambda_2=|\lambda_1|(1-\frac{\lambda_2}{|\lambda_1|})$,
we have
\[
R_{\mathcal F}R_{\mathcal G}
    \le
    \frac{\sqrt{|\mathcal F| |\mathcal G|}(\gbinom{n-1}{k-1}-\sqrt{|\mathcal F| |\mathcal G|})}{\gbinom{n-1}{k-1}(1-\frac{\lambda_2}{|\lambda_1|})}.
\]
\end{proof}

\section{Proof of Theorem~\ref{thm:main}}\label{section-proof-mainth}

Now we are ready to prove Theorem~\ref{thm:main}. Assume to the contrary that
\[
    \delta_d(\mathcal F)\delta_d(\mathcal G)>\gbinom{n-d-1}{k-d-1}^2.
\]
Counting the number of pairs $(S,F)$ satisfying $\left\{(S,F):S\in \gbinom{V}{d},F\in\mathcal{F}, S\le F\right\}$ in two different ways, we obtain
\[
    \sum_{S\in \gbinom{V}{d}}d_S(\mathcal F)
    =
    \gbinom{k}{d}|\mathcal F|.
\]
Thus we have $|\mathcal{F}|\geq \frac{\gbinom{n}{d}}{\gbinom{k}{d}}\delta_d(\mathcal F)$. Similarly,  $|\mathcal{G}|\geq \frac{\gbinom{n}{d}}{\gbinom{k}{d}}\delta_d(\mathcal G)$.
Hence
\[
    X=\sqrt{|\mathcal F||\mathcal G|}
    >\frac{\gbinom{n}{d}}{\gbinom{k}{d}}\gbinom{n-d-1}{k-d-1}:=P.
\]
Thus $X>P$.

We record a useful simplification.

\begin{lemma}\label{lem:P-over-H}
We have
\[
    \frac {P}{\gbinom{n-1}{k-1}}
    =
    \frac{(q^{k-d}-1)(q^n-1)}
         {(q^{n-d}-1)(q^k-1)}.
\]
In particular,
    $0<P<\gbinom{n-1}{k-1}.$
Moreover,
    $1-\frac{P}{\gbinom{n-1}{k-1}}
    =
    \frac{(q^d-1)q^{k-d}(q^{n-k}-1)}
         {(q^{n-d}-1)(q^k-1)}.$
\end{lemma}

\begin{proof}
Using the product formula for Gaussian coefficients and cancelling the common
factors gives 
\begin{equation*}
  \frac {P}{\gbinom{n-1}{k-1}}=\frac{
    \gbinom{n}{d}\gbinom{n-d-1}{k-d-1}
    }{\gbinom{k}{d}\gbinom{n-1}{k-1}}=
   \frac{(q^{k-d}-1)(q^n-1)}
         {(q^{n-d}-1)(q^k-1)}. 
\end{equation*}
Since $n>k$ and $n-d>k-d$, this ratio is strictly less than $1$ and positive.

Finally,
\begin{equation*}
1-\frac{P}{\gbinom{n-1}{k-1}}=1-
    \frac{(q^{k-d}-1)(q^n-1)}
         {(q^{n-d}-1)(q^k-1)}=\frac{
    (q^d-1)q^{k-d}(q^{n-k}-1)
    }{(q^{n-d}-1)(q^k-1)}.   
\end{equation*}
\end{proof}

\subsection{A combinatorial lower bound on $\Omega$}

Define
\[
    c=
    \sqrt{\frac{\delta_d(\mathcal G)}{\delta_d(\mathcal F)}}.
\]
For any $S,T\in \gbinom{V}{d}$, we have
\[
    c\,d_S(\mathcal F)
    \ge
    c\,\delta_d(\mathcal F)
    =
    \sqrt{\delta_d(\mathcal F)\delta_d(\mathcal G)}
    >
    \gbinom{n-d-1}{k-d-1},
\]
and similarly $c^{-1}d_T(\mathcal G)>\gbinom{n-d-1}{k-d-1}$.
Therefore
\[
    \sum_{\substack{S,T\in \gbinom{V}{d}\\S\cap T=\{0\}}}
    \left(c\,d_S(\mathcal F)-\gbinom{n-d-1}{k-d-1}\right)
    \left(c^{-1}d_T(\mathcal G)-\gbinom{n-d-1}{k-d-1}\right)>0.
\]
Expanding the left-hand side of the preceding inequality, we obtain
\[
\begin{aligned}
&\sum_{S\cap T=\{0\}}d_S(\mathcal F)d_T(\mathcal G)-c\gbinom{n-d-1}{k-d-1}\sum_{S\cap T=\{0\}}d_S(\mathcal F)-c^{-1}\gbinom{n-d-1}{k-d-1}\sum_{S\cap T=\{0\}}d_T(\mathcal G)\\
&+\gbinom{n-d-1}{k-d-1}^2|\left\{(S,T):S,T\in \gbinom{V}{d},\ S\cap T=\{0\}\right\}|>0.
\end{aligned}
\]

By definition,
\[
    \Omega
    =
    \sum_{S\cap T=\{0\}}d_S(\mathcal F)d_T(\mathcal G).
\]
For a fixed $S\in\gbinom{V}{d}$, Lemma~\ref{lem:disjoint-count} gives
\[
    |\left\{T\in\gbinom{V}{d}:S\cap T=\{0\}\right\}|
    =q^{d^2}\gbinom{n-d}{d}.
\]
Thus
\[
    \sum_{S\cap T=\{0\}}d_S(\mathcal F)
    =
    q^{d^2}\gbinom{n-d}{d}\sum_{S\in \gbinom{V}{d}}d_S(\mathcal F)
    =q^{d^2}\gbinom{n-d}{d}\gbinom{k}{d}|\mathcal F|,
\] and similarly
$\sum_{S\cap T=\{0\}}d_T(\mathcal G)=q^{d^2}\gbinom{n-d}{d}\gbinom{k}{d}|\mathcal G|$. Moreover,
\[
    |\left\{(S,T):S,T\in \gbinom{V}{d},\ S\cap T=\{0\}\right\}|
    =\gbinom{n}{d} q^{d^2}\gbinom{n-d}{d}.
\]

Therefore
\[
\Omega
    >\gbinom{n-d-1}{k-d-1}q^{d^2}\gbinom{n-d}{d}\gbinom{k}{d}(c|\mathcal F|+c^{-1}|\mathcal{G}|)-\gbinom{n-d-1}{k-d-1}^2 q^{d^2}\gbinom{n-d}{d}\gbinom{n}{d}
\]
By the arithmetic-geometric mean inequality, $c|\mathcal F|+c^{-1}|\mathcal{G}|\ge 2\sqrt{|\mathcal F||\mathcal G|}$. Recall that $X=\sqrt{|\mathcal F||\mathcal G|}$. Therefore, we obtain the following lower bound on $\Omega$.

\begin{corollary}\label{lowerbound-omerga}
Under the conditions stated in Theorem \ref{thm:main} and suppose that
    $\delta_{d}(\mathcal{F})\delta_{d}(\mathcal{G})>\gbinom{n-d-1}{k-d-1}^2$,
we have 
    \[
    \Omega>\alpha_1X-\alpha_2,
\]
where
\[
    \alpha_1
    =
    2\gbinom{n-d-1}{k-d-1}
    q^{d^2}\gbinom{n-d}{d}\gbinom{k}{d},\quad
    \alpha_2
    =
    \gbinom{n-d-1}{k-d-1}^2
    q^{d^2}\gbinom{n-d}{d}\gbinom{n}{d}.
\]
\end{corollary}

\subsection{A spectral upper bound on $\Omega$}

Recall that  
\[
    \chi_{\mathcal F}=a_0+a_1+\cdots+a_k,
    \qquad
    \chi_{\mathcal G}=b_0+b_1+\cdots+b_k,
    \qquad
    a_j,b_j\in U_j,
\]
where $0\leq j\leq k$.
Since $\mathcal F$ and $\mathcal G$ are cross-intersecting, we have $\chi_{\mathcal F}^{\top}A\chi_{\mathcal G}=0$.
It follows that $\sum_{j=0}^k\lambda_j a_j^{\top}b_j=0$. 
By Lemma~\ref{lem:mu-expression}, we have 
$\Omega=\sum_{j=0}^d \mu_j a_j^{\top}b_j$.

We now eliminate the $U_1$-term.  Put $ M:=\frac{\mu_1}{\lambda_1}$.
Since $\lambda_1=-q^{-k}\gbinom{n-k-1}{k-1}<0$
and $\mu_1
    =
    -
    q^{k+(d-1)^2-1}
    \gbinom{k-1}{d-1}
    \gbinom{n-d-1}{d-1}
    \gbinom{n-d-1}{k-d}<0$, 
we have 
\[
    M=\frac{\mu_1}{\lambda_1}
    =
    q^{2k+(d-1)^2-1}
    \frac{
    \gbinom{k-1}{d-1}
    \gbinom{n-d-1}{d-1}
    \gbinom{n-d-1}{k-d}
    }{
    \gbinom{n-k-1}{k-1}
    }>0.
\]
Subtracting $M\sum_{j=0}^k\lambda_j a_j^{\top}b_j=0$
from the expression for $\Omega$, we obtain
\[
    \Omega
    =
    (\mu_0-M\lambda_0)a_0^{\top}b_0
    +
    \sum_{j=2}^k c_j a_j^{\top}b_j, \text{ where }     c_j=
    \begin{cases}
    \mu_j-M\lambda_j, & 2\le j\le d,\\[3pt]
    -M\lambda_j, & d<j\le k.
    \end{cases}
\]
Recall that $X=\sqrt{|\mathcal F||\mathcal G|}$. Since $a_0^{\top}b_0=\frac{X^2}{\gbinom{n}{k}}$, we obtain
\[
    \Omega
    =
    \alpha_3X^2
    +
    \sum_{j=2}^k c_j a_j^{\top}b_j, \text{ where } \alpha_3=\frac{\mu_0-M\lambda_0}{\gbinom{n}{k}}.
\]

We now simplify $\alpha_3$. Recall that $P=\frac{\gbinom{n}{d}}{\gbinom{k}{d}}\gbinom{n-d-1}{k-d-1}$. 
\begin{lemma}\label{lem:alpha3}
\[
    \alpha_3
    =
    \frac{\gbinom{n-d-1}{k-d-1}q^{d^2}\gbinom{n-d}{d}\gbinom{k}{d}}{\gbinom{n-1}{k-1}}
    \left(2-\frac{P}{\gbinom{n-1}{k-1}}\right).
\]
\end{lemma}
\begin{proof}
First, using $\mu_0
    =
    q^{d^2}
    \gbinom{k}{d}
    \gbinom{n-d}{d}
    \gbinom{n-d}{k-d}$ and 
    $\gbinom{n}{d}\gbinom{n-d}{k-d}
    =
    \gbinom{n}{k}\gbinom{k}{d}$,
we get 
\[
\frac{\mu_0}{\gbinom{n}{k}}=\frac{q^{d^2}
    \gbinom{k}{d}
    \gbinom{n-d}{d}
    \gbinom{n-d}{k-d}}{\gbinom{n}{k}}=\frac{q^{d^2}
    \gbinom{n-d}{d}\gbinom{k}{d}^2}{\gbinom{n}{d}}.
\]
Since $P=\frac{\gbinom{n}{d}}{\gbinom{k}{d}}\gbinom{n-d-1}{k-d-1}$, this can be rewritten as 
\[
\frac{\mu_0}{\gbinom{n}{k}}=\frac{\gbinom{n-d-1}{k-d-1}q^{d^2}
    \gbinom{n-d}{d}\gbinom{k}{d}}{P}.
\]

Next we compute $M\lambda_0/\mu_0$.  Since
$\lambda_0=\gbinom{n-k}{k}$, 
we have
\[
\begin{aligned}
    \frac{M\lambda_0}{\mu_0}
    &=
    q^{2(k-d)}
    \frac{
    \gbinom{k-1}{d-1}
    \gbinom{n-d-1}{d-1}
    \gbinom{n-d-1}{k-d}
    \gbinom{n-k}{k}
    }{
    \gbinom{k}{d}
    \gbinom{n-d}{d}
    \gbinom{n-d}{k-d}
    \gbinom{n-k-1}{k-1}
    }                                                        \\
    &=
    \frac{
    q^{2(k-d)}(q^d-1)^2(q^{n-k}-1)^2
    }{
    (q^{n-d}-1)^2(q^k-1)^2
    }.
\end{aligned}
\]
By Lemma~\ref{lem:P-over-H}, $1-\frac{P}{\gbinom{n-1}{k-1}}
    =
    \frac{(q^d-1)q^{k-d}(q^{n-k}-1)}
         {(q^{n-d}-1)(q^k-1)}$. 
Therefore
\[
    \frac{M\lambda_0}{\mu_0}
    =
    \left(1-\frac{P}{\gbinom{n-1}{k-1}}\right)^2.
\]

Consequently,
\[
\begin{aligned}
    \alpha_3
    &=
    \frac{\mu_0-M\lambda_0}{\gbinom{n}{k}}=
    \frac{\mu_0}{\gbinom{n}{k}}
    \left(1-\frac{M\lambda_0}{\mu_0}\right)                               \\
    &=\frac{\gbinom{n-d-1}{k-d-1}q^{d^2}
    \gbinom{n-d}{d}\gbinom{k}{d}}{P}
    \left[
    1-\left(1-\frac{P}{\gbinom{n-1}{k-1}}\right)^2
    \right]                                             \\
    &=\frac{\gbinom{n-d-1}{k-d-1}q^{d^2}
    \gbinom{n-d}{d}\gbinom{k}{d}}{\gbinom{n-1}{k-1}}\left(2-\frac{P}{\gbinom{n-1}{k-1}}\right).
\end{aligned}
\]
\end{proof}

Recall that 
\[
c_j=
\begin{cases}
    \mu_j-M\lambda_j, & 2\le j\le d,\\[3pt]
    -M\lambda_j, & d<j\le k.
\end{cases}
\]
Set
\[
    W=\max_{2\le j\le k}|c_j|.
\]
We need the following upper bound on 
$W$.

\begin{lemma}\label{lem:W}
We have
\[
    W\le M\lambda_2, \text{ where } M=\frac{\mu_1}{\lambda_1}.
\]
More explicitly,
\[
    M\lambda_2
    =
    q^{(d-1)^2}
    \frac{
    \gbinom{k-1}{d-1}
    \gbinom{n-d-1}{d-1}
    \gbinom{n-d-1}{k-d}
    \gbinom{n-k-2}{k-2}
    }{
    \gbinom{n-k-1}{k-1}
    }.
\]
\end{lemma}
\begin{proof}
Suppose that $d=1$. By Lemma \ref{lem:mu-expression}, we have $\mu_j=0$ for every $j\ge 2$.  Hence $c_j=-M\lambda_j$ for all $j\geq 2$.
By Lemma~\ref{lem:lambda-comparison}, $\max_{2\le j\le k}|\lambda_j|=\lambda_2$.
Therefore
\[
    W=M\lambda_2.
\]

Now suppose $d\ge 2$.  Then $c_2=\mu_2-M\lambda_2$. We first show that $M\lambda_2>4\mu_2$. A direct calculation gives
\[
\begin{aligned}
    \frac{M\lambda_2}{\mu_2}
    &=q^{-2(k-d)}
    \left(
    \frac{(q^{k-1}-1)(q^{n-d-1}-1)}
         {(q^{d-1}-1)(q^{n-k-1}-1)}
    \right)^2  \\
    &>q^{-2(k-d)}q^{4(k-d)}
    =
    q^{2(k-d)}
    \ge 4,
\end{aligned}
\]
where the first inequality holds as $k>d$ and 
$n-d-1>n-k-1$. Thus
\[
    |c_2|=M\lambda_2-\mu_2>0.
\]

We show that the sequence $(|\mu_j|)_{j\geq 2}$ is decreasing.
For $2\le j<d\leq k$, a direct calculation gives
\[
\begin{aligned}
    \frac{|\mu_{j+1}|}{|\mu_j|}
    &=
    q^{k-2d+j}
    \frac{q^{d-j}-1}{q^{k-j}-1}
    \frac{q^{d-j}-1}{q^{n-d-j}-1}
    \frac{q^{n-k-j}-1}{q^{n-d-j}-1}<q^{k-2d+j}q^{-(k-d)}
    =
    q^{j-d}
    <1.
\end{aligned}
\]
Therefore
\[
    |\mu_j|\le \mu_2,
    \quad
    \text{for } j\ge 2.
\]

Now let $j\ge 3$. If $j\le d$, since \(\mu_j\) and \(\lambda_j\) have the same sign for every
\(j\le d\), and hence 
\[
\begin{aligned}
    |c_j|
    =
    |\mu_j-M\lambda_j|&=\big||\mu_j|-M|\lambda_j|\big|    
    \le
    \max\{|\mu_j|,M|\lambda_j|\}                      \\
    &\le\max\{|\mu_2|,M|\lambda_2|\}=M\lambda_2.
\end{aligned}
\]
If $j>d$, then 
\[
    |c_j|=M|\lambda_j|
    \le
    M\lambda_2.
\]
This completes the proof.
\end{proof}
   
We now use $W$ to bound the spectral error term. Combining Lemmas \ref{lem:alpha3} and \ref{lem:W}, we obtain the following upper bound on $\Omega$.

\begin{corollary}\label{lem:spectral-upper}
Under the conditions stated in Theorem \ref{thm:main}, we have
\[
    \Omega
    \le
    \alpha_3X^2
    +
    \frac{W}{\gbinom{n-1}{k-1}(1-\frac{\lambda_2}{|\lambda_1|})}X\left(\gbinom{n-1}{k-1}-X\right),
\]
where
\[
\alpha_3
    =
    \frac{\gbinom{n-d-1}{k-d-1}q^{d^2}\gbinom{n-d}{d}\gbinom{k}{d}}{\gbinom{n-1}{k-1}}
\left(2-\frac{P}{\gbinom{n-1}{k-1}}\right).
\]
\end{corollary}

\begin{proof}
From the spectral reduction, $\Omega
    =
    \alpha_3X^2
    +
    \sum_{j=2}^k c_j a_j^{\top}b_j$.
Therefore
\[
\begin{aligned}
    \left|
    \sum_{j=2}^k c_j a_j^{\top}b_j
    \right|
    &\le
    W\sum_{j=2}^k |a_j^{\top}b_j|                         \\
    &\le
    W
    \left(\sum_{j=2}^k\|a_j\|^2\right)^{1/2}
    \left(\sum_{j=2}^k\|b_j\|^2\right)^{1/2}   =
    W R_{\mathcal F}R_{\mathcal G}.
\end{aligned}
\]
By Lemma~\ref{lem:cross-stability}, $R_{\mathcal F}R_{\mathcal G}
    \le
    \frac{X(\gbinom{n-1}{k-1}-X)}{\gbinom{n-1}{k-1}(1-\frac{\lambda_2}{|\lambda_1|})}$.
Hence
\[
    \Omega
    \le
    \alpha_3X^2
    +
    \frac{W}{\gbinom{n-1}{k-1}(1-\frac{\lambda_2}{|\lambda_1|})}X\left(\gbinom{n-1}{k-1}-X\right). 
\]
\end{proof}

\subsection{Finishing the proof}

Put $K:=\frac{W}{\gbinom{n-1}{k-1}(1-\frac{\lambda_2}{|\lambda_1|})}$. We next derive a convenient upper bound on $K$.

\begin{lemma}\label{lem:K-comparison}
We have 
\[
    K<
    \frac{\gbinom{n-d-1}{k-d-1}q^{d^2}
    \gbinom{n-d}{d}\gbinom{k}{d}}{\gbinom{n-1}{k-1}}\left(1-\frac{P}{\gbinom{n-1}{k-1}}\right).
\]
\end{lemma}
\begin{proof}
By Lemma~\ref{lem:W},
\[
    W\le M\lambda_2.
\]    
We compute
\[
\begin{aligned}
    \frac{M\lambda_2}{\gbinom{n-d-1}{k-d-1}q^{d^2}
    \gbinom{n-d}{d}\gbinom{k}{d}}&=\frac{
    q^{(d-1)^2}
    \frac{
    \gbinom{k-1}{d-1}
    \gbinom{n-d-1}{d-1}
    \gbinom{n-d-1}{k-d}
    \gbinom{n-k-2}{k-2}
    }{
    \gbinom{n-k-1}{k-1}
    }
    }{
    \gbinom{n-d-1}{k-d-1}
    q^{d^2}\gbinom{n-d}{d}\gbinom{k}{d}
    }                 \\
     &=
    q^{-2d+1}
    \frac{\gbinom{k-1}{d-1}}{\gbinom{k}{d}}
    \frac{\gbinom{n-d-1}{d-1}}{\gbinom{n-d}{d}}
    \frac{\gbinom{n-d-1}{k-d}}{\gbinom{n-d-1}{k-d-1}}
    \frac{\gbinom{n-k-2}{k-2}}{\gbinom{n-k-1}{k-1}}              \\
    &=
    q^{-2d+1}
    \frac{q^d-1}{q^k-1}
    \frac{q^d-1}{q^{n-d}-1}
    \frac{q^{n-k}-1}{q^{k-d}-1}
    \frac{q^{k-1}-1}{q^{n-k-1}-1}.
\end{aligned}
\]
By Lemma~\ref{lem:P-over-H},
\[
    1-\frac{P}{\gbinom{n-1}{k-1}}
    =
    \frac{(q^d-1)q^{k-d}(q^{n-k}-1)}
         {(q^{n-d}-1)(q^k-1)}.
\]
Therefore
\[
\frac{M\lambda_2}{\gbinom{n-d-1}{k-d-1}q^{d^2}
    \gbinom{n-d}{d}\gbinom{k}{d}\left(1-\frac{P}{\gbinom{n-1}{k-1}}\right)}=q^{1-k-d}
    \frac{(q^d-1)(q^{k-1}-1)}
         {(q^{k-d}-1)(q^{n-k-1}-1)}.
\]
By Lemma \ref{lem:lambda-comparison}, we recall that 
\[
r=\frac{\lambda_2}{|\lambda_1|}
    =
    q^{1-k}
    \frac{q^{k-1}-1}{q^{n-k-1}-1}.
\]
Hence 
\[
\frac{M\lambda_2}{\gbinom{n-d-1}{k-d-1}q^{d^2}
    \gbinom{n-d}{d}\gbinom{k}{d}\left(1-\frac{P}{\gbinom{n-1}{k-1}}\right)}=r\cdot
    q^{-d}\frac{q^d-1}{q^{k-d}-1}<r.
\]
Since $W\le M\lambda_2$ and $r=\frac{\lambda_2}{|\lambda_1|}$, it follows that
\[
\frac{W}{(1-r)\gbinom{n-d-1}{k-d-1}q^{d^2}
    \gbinom{n-d}{d}\gbinom{k}{d}\left(1-\frac{P}{\gbinom{n-1}{k-1}}\right)}<\frac{r}{1-r},
\]
which is less than or equal to $1$ as $r\leq 1/2$ by Lemma~\ref{lem:lambda-comparison}. This completes the proof.
\end{proof}

We now complete the proof of Theorem~\ref{thm:main}. Recall that $X=\sqrt{|\mathcal{F}||\mathcal{G}|}$.
Combining Corollaries \ref{lowerbound-omerga} and \ref{lem:spectral-upper}, we obtain
\[
0<
    \alpha_3X^2-\alpha_1X+\alpha_2+KX\left(\gbinom{n-1}{k-1}-X\right).
\]
We now factor the quadratic polynomial $\alpha_3t^2-\alpha_1t+\alpha_2$. Recall that
\[
\alpha_1
    =
    2\gbinom{n-d-1}{k-d-1}
    q^{d^2}\gbinom{n-d}{d}\gbinom{k}{d},\quad
    \alpha_2
    =
    \gbinom{n-d-1}{k-d-1}^2
    q^{d^2}\gbinom{n-d}{d}\gbinom{n}{d}, 
\]
and
\[
\alpha_3
    =
    \frac{\gbinom{n-d-1}{k-d-1}q^{d^2}\gbinom{n-d}{d}\gbinom{k}{d}}{\gbinom{n-1}{k-1}}
\left(2-\frac{P}{\gbinom{n-1}{k-1}}\right).
\]
Since $\alpha_3 \gbinom{n-1}{k-1}^2-\alpha_1\gbinom{n-1}{k-1}+\alpha_2=\gbinom{n-d-1}{k-d-1}q^{d^2}\gbinom{n-d}{d}\gbinom{k}{d}\left(\left(2-\frac{P}{\gbinom{n-1}{k-1}}\right)\gbinom{n-1}{k-1}-2\gbinom{n-1}{k-1}+\frac{\gbinom{n}{d}\gbinom{n-d-1}{k-d-1}}{\gbinom{k}{d}}\right)$ and $P=\frac{\gbinom{n}{d}\gbinom{n-d-1}{k-d-1}}{\gbinom{k}{d}}$,
it follows that $\gbinom{n-1}{k-1}$ is a root of the quadratic polynomial $\alpha_3t^2-\alpha_1t+\alpha_2$. The other root is 
\[
s=\frac{\alpha_2}{\alpha_3\gbinom{n-1}{k-1}}=\frac{P}{2-\frac{P}{\gbinom{n-1}{k-1}}}.
\]
By Lemma \ref{lem:P-over-H}, $0<\frac{P}{\gbinom{n-1}{k-1}}<1$. Then we have $0<s<P<\gbinom{n-1}{k-1}$. Moreover, 
\[
\frac{\gbinom{n-d-1}{k-d-1}q^{d^2}
    \gbinom{n-d}{d}\gbinom{k}{d}}{\gbinom{n-1}{k-1}}\left(1-\frac{P}{\gbinom{n-1}{k-1}}\right)=\frac{\alpha_3}{2-\frac{P}{\gbinom{n-1}{k-1}}}\left(1-\frac{P}{\gbinom{n-1}{k-1}}\right)
    =\alpha_3\frac{P-s}{P}.
\]
Thus the inequality in Lemma \ref{lem:K-comparison} can be written as $K< \alpha_3\frac{P-s}{P}$.

By Lemma \ref{lem:cross-stability}, we have the Tokushige bound $X\leq \gbinom{n-1}{k-1}$. Then 
\[
\begin{aligned}
0<&\alpha_3X^2-\alpha_1X+\alpha_2+KX\left(\gbinom{n-1}{k-1}-X\right) \\
\leq &\alpha_3\left(X-\gbinom{n-1}{k-1}\right)(X-s)+\alpha_3\frac{P-s}{P}X\left(\gbinom{n-1}{k-1}-X\right)\\
=&\alpha_3\left(\gbinom{n-1}{k-1}-X\right)\left(\frac{P-s}{P}X-X+s\right)\\
=&\alpha_3\left(\gbinom{n-1}{k-1}-X\right)\frac{s(P-X)}{P}
\end{aligned}
\]
The assumption $\delta_d(\mathcal F)\delta_d(\mathcal G)>\gbinom{n-d-1}{k-d-1}^2$ now implies that $X>P$. Since $s<P<X$, we obtain
\[
\alpha_3\left(\gbinom{n-1}{k-1}-X\right)\left(s-\frac{s}{P}X\right)\leq 0.
\]
This contradiction shows that
    $\delta_d(\mathcal F)\delta_d(\mathcal G)>\gbinom{n-d-1}{k-d-1}^2$
is impossible.  Hence
\[
    \delta_d(\mathcal F)\delta_d(\mathcal G)
    \le
    \gbinom{n-d-1}{k-d-1}^2.
\]
The proof of Theorem~\ref{thm:main} is now complete.

\end{document}